\documentclass[10pt]{article}
\usepackage{amsfonts}
\usepackage{amssymb}
\newcommand{\text}{\mbox}
\newcommand{\oc}{{\sf{co}}}
\newcommand{\Id}{\text{Id}}
\newcommand{\C}{{\bf C}}

\renewcommand{\H}{{\bf H}}
\newcommand{\M}{{\bf M}}
\newcommand{\1}{{\bf 1}}

\newcommand{\ran}{\text{ran}}
\newcommand{\Comp}{{\text{Comp}}}
\begin{document}

\begin{center}{\large \sc Topological Orthoalgebras}\\

Alexander Wilce\\ Department of Mathematical Sciences\\
Susquehanna University\\Selinsgrove, PA 17870 \ \ e-mail: wilce@susqu.edu\\ \end{center}

\begin{abstract}  We define topological orthoalgebras
(TOAs), and study their properties. While every topological orthomodular lattice
is a TOA, the lattice of projections of a Hilbert space (in its norm or strong operator
topology) is an example of a lattice-ordered TOA that is not a topological lattice.
On the other hand, we show that every compact Boolean TOA is a
topological Boolean algebra.  We also show that a compact TOA
in which $0$ is an isolated point is atomic and of finite height. We identify and study a
particularly tractable class of TOAs, which we call {\em stably ordered}: those in which the upper-set generated by an open set is open. This includes all topological OMLs, and also the projection lattices of Hilbert spaces. Finally, we obtain a version of the Foulis-Randall representation theory for stably-ordered topological orthoalgebras. Keywords: Topological orthoalgebra, topological orthomodular lattice. AMS subject classification: 03G12, 06C15, 06F30, 54H12. 
\end{abstract}

\underline{\bf 0. Introduction}

An orthomodular poset (OMP) is an orthocomplemented poset 
$(L,\leq,')$ in which, for all $a \leq b \in L$, 
\begin{equation}
\text{$b \wedge a'$ exists} \end{equation}
(whence also $b' \vee a = (b \wedge a')'$ exists), and
\begin{equation}
(b \wedge a') \vee a = b.\end{equation}
An {\em orthomodular lattice} (OML) is a lattice-ordered OMP.
Evidently, the ``orthomodular identity" (2) is a weak form of distributivity,
and thus every Boolean algebra is an OML. The most important non-Boolean
example is the geometry of closed subspaces of a Hilbert space, or more
generally, of a von Neumann algebra. Orthomodular lattices have been extensively
studied (see, e.g., [9]), both in connection with quantum logic and also as objects of
intrinsic interest. In particular,
there is a fairly developed theory of topological OMLs due to Choe and
Greechie [2,3], Pulmannov\'{a} and Riecanova [13], and others. However, this
does not accommodate projection lattices (in the norm topology),
since the meet and join operations on such a lattice are not continuous. 

In any orthoposet, we may define a partial binary
operation by setting $a \oplus b = a \vee b$ provided that
$a \leq b'$ and the join exists. In the case of the projection 
lattice of a Hilbert space, $a \leq b'$ iff the range of $a$ is orthogonal to that of $b$, in which case $a \oplus b$ is exactly the projection onto the closed linear span of $a$ and $b$. Thus, every orthomoposet $(L,\leq)$ gives rise
to a partial semigroup $(L,\oplus)$. If $L$ is orthomodular, this satisfies the axioms for
an {\em orthoalgebra} [4], given below. Conversely, every orthoalgebra carries a natural
partial ordering, and it can be shown that OMLs are equivalent to
lattice-ordered orthoalgebras. In this paper, we develop a theory of topological orthoalgebras. This subsumes the theory of topological OMLs, in that every topological OML is also a
topological orthoalgebra. On the other hand, the projection
lattice of a Hilbert space, in the operator norm topology, is also
a topological orthoalgebra. 

This paper is in part a sequel to [16]; however, we have included
enough detail to make it reasonably self-contained. \\

\underline{\bf 1. Background on Orthoalgebras and Test Spaces}

If $(L,\leq,0,1)$ is any orthocomplemented poset, we call elements $a$ and $b$ of $L$
{\em orthogonal}, writing $a \perp b$, iff $a \leq b'$. Let us agree in this case to write $a \oplus b$ for the join of $a$ and $b$, whenever the latter exists. This partial operation is commutative and associative in the strong sense that if $a \perp b$ and $(a \oplus b) \perp c$, then  $b \perp c$, $a \perp (b \oplus c)$,  and
$(a \oplus b) \oplus c = a \oplus (b \oplus c)$.  Note also that $a \oplus a' = \1$ for 
every $a \in L$. The following is easily proved:

{\bf 1.1 Lemma:} {\em Let $(L, \leq, ')$ be an orthoposet in which $a \oplus b = a \vee b$ exists whenever $a \leq b'$. Then the following are equivalent: 
\begin{itemize}
\item[(a)] $(L,\leq,')$ is an orthomodular poset,
\item[(b)] $\forall a \in L$, $a \oplus b = 1 \Rightarrow b = a'$, 
\item[(c)] $(L,\oplus)$ is cancellative, i.e., $a \oplus b = a \oplus c \Rightarrow b = c$.
\end{itemize}}

This suggests the following 

{\bf 1.2 Definition:} An {\bf orthoalgebra} is a structure
$(L,\oplus,0,1)$ consisting of a set $L$, a partially-defined
associative, commutative\footnote{Commutativity is also understood in
the strong sense that if $a \oplus b$ is defined, then so is $b \oplus a$,
and the two are equal.} binary operation $\oplus$ on $L$ and two distinguished elements
$0$ and $1$, such that, for all $a \in L$,
\begin{itemize}
\item[(a)] $a \oplus 0 = 0 \oplus a = a$;
\item[(b)] there exists a unique element $a' \in L$ with $a \oplus a' = 1$;
\item[(c)] $a \oplus a$ exists only if $a = 0$.
\end{itemize}
Orthoalgebras were introduced in the early 1980s by D. J. Foulis and C. H. Randall [6] in connection with the problem of defining tensor products of quantum logics.\footnote{The axioms given here is due to A. Golfin.} Further information can be found in [4] and [15].
 For later reference, we mention that an {\em effect algebra} (see, e.g., [1]) is a structure $(L,\oplus,0,1)$ satisfying conditions (a), (b) and (c), but not necessarily (d).

It is straightforward to check that the structure $(L,\oplus,0,1)$ obtained
as above from an orthomodular poset $(L,\leq, \, ' , 0, 1)$ is an orthoalgebra.
Conversely, any orthoalgebra $(L,\oplus,0,1)$ can be partially ordered by
setting $a \leq b$ iff there exists $c \in L$ with $b = a \oplus c$. The
operation $a \mapsto a'$ is an orthocomplementation with respect to this ordering.
Thus, any orthoalgebra gives rise to an orthoposet. Moreover, for
any $a \leq b$ in $L$, there is a unique element $c \in L$ -- namely,
$(b \oplus a')'$ -- such that $b = a \oplus c$. It is usual to call this
element $b \ominus a$. If $L$ is an OMP, this is exactly $b \wedge a'$. In 
this language, 
the orthomodular law (2) becomes 
\begin{equation}
a \leq b \ \Rightarrow \ b = (b \ominus a) \oplus a, \end{equation}
which holds in any orthoalgebra.

{\bf 1.3 Orthomodular Posets as Orthoalgebras.} In general, $a \oplus b$ is not the join, but only a {\em minimal} upper bound, for orthogonal elements $a$ and $b$ of an orthoalgebra $L$. Indeed, one can show that the orthoposet $(L,\leq,',0,1)$ obtained from $(L,\oplus,0,1)$ is an OMP if and only if $a \oplus b = a \vee b$ for all $a, b \in L$. This in turn is equivalent to the condition that if $a, b, c \in L$ with $a \oplus b$, $b \oplus c$ and $c \oplus a$ all defined,
then $a \oplus (b \oplus c)$ is defined [4,15]. This last condition is called {\em orthocoherence} in the literature. Thus, OMPs are essentially the same things as orthocoherent orthoalgebras.  

{\bf 1.4 Intervals and Products.} If $L$ is an orthoalgebra and $a \in L$, then the order-interval
\[[0,a] \ := \ \{ b \in L | 0 \leq b \leq a\}\] 
can be made into an orthoalgebra as follows: define $x \perp_a y$ to mean that $x \perp y$ in $L$ and $x \oplus y \leq a$; in this case, set $x \oplus_a y = x \oplus y$. Then $([0,a],\oplus_a,0,a)$ is an orthoalgebra. If $L_1$ and $L_2$ are any two orthoalgebras (or effect algebras), we can give $L_1 \times L_2$ the structure of an orthoalgebra by setting $(x,y) \perp (u,v)$ iff $x \perp u$ and $y \perp v$, in which case $(x,u) \oplus (y,v) = (x \oplus u,y \oplus v)$. The zero and unit elements in $L_1 \times L_2$ are then $(0,0)$ and $(1,1)$, respectively, and $(x,y)' = (x',y')$. Let $a = (1,0) \in L$: then 
$L_1 \simeq [0,a]$ and $L_2 \simeq [0,a']$ under the mappings $x \mapsto (x,0)$ and $y \mapsto (0,y)$, respectively.

{\bf 1.5 Compatibility.} An orthoalgebra $(L,\oplus,0,1)$ is said to be {\em Boolean} iff the corresponding orthoposet $(L,\leq,',0,1)$ is a Boolean lattice.
Elements $a$ and $b$ of an orthoalgebra $L$ are said to be {\em compatible} (or to {\em commute} iff they generate a Boolean sub-orthoalgebra of $L$. This is
the case iff there exist elements $a_1, b_1$ and $c$ with $a = a_1 \oplus c$,
$b = c \oplus b_1$, and $a \perp b_1$, so that $a_1 \oplus c \oplus b_1$
exists [4]. Equivalently, $a$ and $b$ are compatible iff there exists an element $c \leq a, b$ with $a \perp (b \ominus c)$. The triple $(a_1, c, b_1) = (a \ominus c, c, b \ominus c)$ is then called a {\em Mackey decomposition} for $a$ and $b$. 
If $L$ is Boolean, then every pair of elements $a, b \in L$ has a {\em unique} Mackey decomposition, namely, $(a \ominus b, a \wedge b, b \ominus a)$.  More generally, if $L$ is orthocoherent, i.e., an OMP, then Mackey decompositions are unique where they exist, and again $c = a \wedge b$. 

Let us call elements $a$ and $b$ of an OA $L$ {\em uniquely compatible} iff they have a unique Mackey decomposition. For any $a \in L$, there is a natural mapping $[0,a] \times [0,a'] \rightarrow L$ given by $(x,y) \mapsto x \oplus y$. If this mapping is in fact an isomorphism, $a$ is said to be {\em central}. The {\em center} of $L$ is the set $\C(L)$ of all central elements of $L$. It is not difficult to show that $a$ is central iff $a$ is uniquely compatible with every $b \in L$; thus if $L$ is an OMP, $\C(L)$ is the set of all elements compatible with every element. It can be shown [7] that $\C(L)$ is a Boolean sub-orthoalgebra algebra of $L$. In particular, $L$ is Boolean iff $L = \C(L)$.

We shall call an orthoalgebra {\em simple} iff its center is the trivial Boolean algebra $\{0,1\}$. The projection lattice $L(\H)$ of a Hilbert space is simple. More generally, the lattice of projections of a von Neumann algebra is simple iff the latter is a factor. It should be mentioned that there exist non-orthocoherent simple orthoalgebras in which every pair of elements is compatible [4, Example 3.5]. Note that if $a \in L$ is a central atom, i.e., an atomic element of $\C(L)$, then $[0,a]$ is a simple orthoalgebra.

{\bf 1.6 Test Spaces.} A {\em test space} (see [4] or [15], and references therein) is a pair $(X,{\frak A})$ where $X$ is a non-empty set and ${\frak A}$ is a covering of $X$  by non-empty subsets, called {\em tests}. The usual interpretation is that each set $E \in {\frak A}$ represents the set of  possible {\em outcomes} of some experiment, decision, or physical process; accordingly, any subset of a test is called an {\em event}. The collection of all events of a test space $(X,{\frak A})$ is denoted by ${\cal E}(X,{\frak A})$. Events $A, B \in {\cal E}(X,{\frak A})$  are said to be {\em orthogonal}, written $A \perp B$, if they are disjoint and their union is again an event, and {\em complementary} if they partition a test. $A$ and $B$ are {\em perspective}, written $A \sim B$, if they share a common complementary event. Note that any tests $E, F \in {\frak A}$ are 
perspective (both being complementary to the empty event). 

A test space $(X,{\frak A})$ is {\em algebraic} iff perspective events have {\em exactly the same} complementary events. In this case $\sim$ is an equivalence relation on ${\cal E}$; furthermore, if  $A \perp B$ and $B \sim C$, then $A \perp C$ and $A \cup B \sim A \cup C$. Let $\Pi(X,{\frak A}) = {\cal E}(X,{\frak A})/\sim$, and write $p(A)$ for the $\sim$-equivalence class of an event $A \in {\cal E}(X)$. Then $\Pi$ carries a well-defined partial binary operation $p(A) \oplus p(B) := p(A \cup B)$, defined for orthogonal  pairs of events. Let $0 :=p(\emptyset)$ and $1 := p(E)$, $E \in {\frak A}$: then $(\Pi,\oplus, 0, 1)$ is an orthoalgebra, called the {\em logic} of $(X,{\frak A})$.  

By way of illustration, let $S$ be the unit sphere of a Hilbert space $\H$, and let ${\frak F}$ be the collection of all {\em frames} (maximal orthonormal subsets) of $S$. Then $(S,{\frak F})$ is a test space, the events for which are just pairwise orthogonal sets of unit vectors. Two events are orthogonal iff they span orthogonal subspaces, complementary iff they span complementary subspaces, and perspective iff they span the {\em same} subspace of $\H$. Evidently, $(S,{\frak F})$ is algebraic, with logic isomorphic to $L(\H)$. 

In fact, any orthoalgebra $L$ is canonically representable as the logic of an algebraic test space. If $A = \{a_1,...,a_n\}$ is a finite subset of $L$, let $\bigoplus A = a_1 \oplus \cdots \oplus a_n$, provided this exists. Let $X_L = L \setminus \{0\}$, and let ${\frak A}_L$ denote the set of finite subsets $E \subseteq L$ for which $\bigoplus E$ exists and equals $1$. Then $(X_L, {\frak A}_L)$ is an algebraic test space, the events of which are simply those finite subsets $A$ of $X_L$  for which $\bigoplus A$ exists. The logic of $(X_L, {\frak A}_L)$ is isomorphic to $L$ under  the (well-defined) mapping $p(A) \mapsto \bigoplus A$.\\

\underline{\bf 2.  Topological Orthoalgebras}

{\bf 2.1 Definition:} A {\sl topological orthoalgebra} (TOA) is an
orthoalgebra $(L,\oplus,0,1)$ equipped with a topology making the relation
$\perp \subseteq L \times L$ closed, and the mappings $\oplus : \perp
\rightarrow L$ and $' : L \rightarrow L$, continuous.

{\bf Remark:} One could define a topological effect algebra in just the same way. We shall not pursue this further, except to note that the 
following would apply verbatim to that context. 

{\bf 2.2 Lemma:} {\em Let $(L,\oplus,0,1)$ be a topological
orthoalgebra. Then
\begin{itemize}
\item[(a)] The order relation $\leq$ is closed in $L \times L$
\item[(b)] $L$ is a Hausdorff space.
\item[(c)] The mapping $\ominus : \leq \rightarrow L$ is
continuous.                                           \end{itemize}}

{\bf Proof}: For (a), notice that $a \leq b$ iff $a \perp b'$. Thus, $\leq = f^{-
1}(\perp)$ where $f : L \times L \rightarrow L \times L$ is the continuous
mapping $f(a,b) = (a,b')$. Since $\perp$ is closed, so is $\leq$. That $L$ is
Hausdorff now follows by standard arguments (cf. [8, Chapter VII] or [11]). Finally, since
$b \ominus a = (b \oplus a')'$, and $\oplus$ and $'$ are both continuous,
$\ominus$ is also continuous. $\Box$

Topologically, TOAs are much freer objects than TOMLs. In fact, any Hausdorff space can be embedded in a TOA by the following construction:

{\bf 2.3 Example:} Let $X$ be any Hausdorff space and let $Y$ be a (disjoint) homeomorphic copy of $X$.  Let $f : X \rightarrow Y$ be a homeomorphism. Select a point $0 \in X$, and set $1 := f(0) \in Y$. Let $L = X \cup Y$ (with the obvious topology). For $x, y \in L$, set $x \perp y$ iff (i) $(x,y)$ or $(y,x)$ belongs to the graph of $f$, or (ii)  $x$ or $y$ is $0$. In the former case, set $x \oplus y = 1$; in the latter cases, set $0 \oplus y = y$ and $x \oplus 0 = x$. It is easily checked that $L$ is then a topological orthoalgebra, with orthocomplementation given by $x' = f(x)$ for $x \in X$ and $y' = f^{-1}(y)$ for $y \in Y$. Note that order-theoretically, this is just the horizontal sum of the Boolean algebras $\{0,x,f(x),1\}$, indexed by $x \in X$.

It is straightforward to show that any Cartesian product of discrete
orthoalgebras, with the product topology, is a compact TOA. Also, any topological OML is an example of a topological OA, since in that setting
$a \perp b$ iff  $a \leq b'$ iff $a = a\wedge b'$, which is obviously a
closed relation when $L$ is Hausdorff and $\wedge, '$ are continuous.
However, there are simple and important examples of lattice-ordered TOAs that are not
TOMLs.

{\bf 2.4 Example:} $L$ be the horizontal sum of four-element Boolean algebas
$L_{x} = \{0,x,x',1\}$ with $x$ (and hence, $x'$) parametrized by a non-degenerate real
interval $[a,b]$. Topologize this as two disjoint copies, $I$ and $I'$, of $[a,b]$ plus
two isolated points $0$ and $1$: Then the orthogonality relation is obviously
closed, and $\oplus$ is obviously continuous; however, if we let $x
\rightarrow x_o$ (with $x \not = x_{o})$ in $I$, then we have $x \wedge x_{o}
= 0$ yet $x_o \wedge x_{o} = x_{o}$; hence, $\wedge$ is not continuous. (Note that
this is a special case of Example 2.3.)

{\bf 2.5 Example:} Let $\H$ be a Hilbert space, and let $L = L(\H)$ be the
space of projection operators on $\H$, with the operator-norm topology. The relation
Since multiplication is continuous, the relation $P \perp Q$ iff $PQ = QP = 0$ is closed;
since addition and subtraction are continuous, the partial operation $P, Q \mapsto P \oplus Q := P+Q$ is
continuous on $\perp$, as is the operation $P \mapsto P' := \1 - P$. So $L(\H)$ is a lattice-ordered
topological orthoalgebra. It is not, however, a topological lattice. Indeed, if $Q$ is a non-trivial
 projection, choose unit vectors $x_n$ not
lying in $\ran(Q)$ that converge to a unit vector in $x \in \ran(Q)$. If $P_n$ is the projection generated by $x_n$ and $P$, that generated by with $x$, then $P_n \rightarrow P$. But
$P_n \wedge Q = 0$, while $P \wedge Q = P$.
\footnote{It is worth remarking that $L(\H)$ {\em is} a  (non-Hausdorff) topological OML in its order topology [14].}  

Note that if $\H$ is finite-dimensional, then $L(\H)$ is closed,
and hence, compact in the unit sphere of ${\cal B}(\H)$.  It is a basic result of the theory of topological orthomodular lattices that a compact TOML is totally disconnected [2, Lemma 3]. In strong contrast to this, the subspace of $L(\H)$
consisting of projections of a given rank is a manifold.

For the balance of this section, we concentrate on compact TOAs. We begin with a simple 
completeness result, of which we make no use in the sequel. An ordered set $L$ is {\em order-complete} iff arbitrary {\em directed} joins and meets 
exist in $L$. It is a standard fact that any compact ordered space with a closed order 
is order-complete. (Indeed, any such ordered space is isomorphic to a closed subspace of 
a cartesian power of $[0,1]$ in its product order and topology [12, Corollary VII.1.3]) A subset $J$ of an orthoalgebra $L$ is said to be {\em jointly orthogonal} iff every finite set $A \subseteq J$ is pairwise-orthogonal and has a well-defined orthogonal sum 
$\bigoplus A = a_1 \oplus \cdots \oplus a_n$. We say that $J$ is {\em summable} iff 
the directed $\{\bigoplus A | A \subseteq J\}$ of all finite partial sums of elements of $J$ has a join in $L$, in which case, we write $\bigoplus J$ for this sum. $L$ is {\em complete} iff every jointly-orthogonal subset of $L$ is summable.  Since the set of finite partial sums of a jointly orthogonal set is clearly directed, it follows that an 
order-complete orthoalgebra is a complete orthoalgebra. In particular, then, we have the 

{\bf 2.6 Proposition:} {\em Any compact topological OA is a complete orthoalgebra.}

If $L$ is any orthoalgebra, let 
\[\M(L) := \{ (a,c,b) \in L^{3} | c \leq a, \ c \leq b, \ \text{and} \ a \ominus c \perp
b\}.\] 
In other words, $(a,c,b) \in \M(L)$ iff $(a \ominus c, c, b \ominus c)$ is a Mackey decomposition for $a$ and $b$.

{\bf 2.7 Lemma:} {\em $\M(L)$ is closed in $L^{3}$.}

{\bf Proof}: Just note that
$\M(L) = (\geq \times L) \cap (L \times \leq) \cap (\ominus \times
\Id)^{-1}(\perp)$.
Since $\leq$ is closed and $\ominus : \leq \rightarrow L$ is continuous,
this is also closed. $\Box$

{\bf 2.8 Proposition:} {\em A compact Boolean topological orthoalgebra is a
topological lattice, and hence, a compact topological Boolean algebra.}

{\bf Proof}: If $L$ is Boolean, then $\M(L)$ is, up to a permutation, the
graph of the mapping $a,b \mapsto a \wedge b$. Thus, by Lemma 2.6,
$\wedge$ has a closed graph. Since $L$ is
compact, this suffices to show that $\wedge$ is continuous.\footnote{Recall
here that if $X$ and $Y$ are compact and the graph $G_f$ of $f : X \rightarrow Y$ is
closed, then $f$ is continuous. Indeed, let $F \subseteq Y$ be closed. Then
$f^{-1}(F) = \pi_{1}((X \times F) \cap G_{f})$, where $\pi_{1}$ is projection
on the first factor. Since $X$ and $Y$ are compact, $\pi_{1}$ sends closed
sets to closed sets.} It now follows from the continuity of $'$ that $\vee$
is also continuous. $\Box$

For any orthoalgebra $L$, let $\Comp(L)$ be the set of all
compatible pairs in $L$, and for any fixed $a \in L$, let $\Comp(a)$ be the set
of elements compatible with $a$. 

{\bf 2.9 Proposition:} {\em Let $L$ be a compact TOA. Then
\begin{itemize}
\item[(a)] $\Comp(L)$ is closed;
\item[(b)] For every $b \in L$, $\Comp(b)$ is closed;
\item[(c)] The closure of a pairwise compatible set is pairwise compatible;
\item[(d)] A maximal pairwise compatible set is closed.
\end{itemize}}

{\bf Proof}: $\Comp(L) = (\pi_{1} \times \pi_{3})(\M(L))$. Since $\M(L)$ is closed,
and hence compact, and $\pi_{1} \times \pi_{3}$ is continuous, $\Comp(L)$ is
also compact, hence closed. For (b), note that $\Comp(b) = \pi_{1}(\Comp(L) \cap
(L \times \{b\}))$. Since $\Comp(L)$ is closed, so is $\Comp(L) \cap (L \times
\{b\})$; hence, its image under $\pi_{1}$ is also closed (remembering
here that $L$ is compact). For (c), suppose $M \subseteq L$ is pairwise
compatible. Then $M \times M \subset \Comp(L)$. By part (a), $\Comp(L)$ is closed,
so we have
\[\overline{M} \times \overline{M} \subseteq \overline{M \times M} \subseteq
\Comp(L),\]
whence, $M$ is again pairwise compatible. Finally, for (d), if $M$ is
a maximal pairwise compatible set, then the fact that $M \subseteq
\overline{M}$ and $\overline{M}$ is also pairwise compatible entails that
$M = \overline{M}$.
$\Box$

As mentioned in section 1, there exist non-orthocoherent orthoalgebras for which $\Comp(L) = L \times L$, i.e., every pair of elements are compatible ]. However, in an OMP, pairwise compatible elements always generate a Boolean sub-OML. In this setting, $\Comp(L) = \C(L)$, the center of $L$; thus, we recover from 2.9 the fact (not hard to prove directly; see  [2]) that the center of a compact TOML is a compact Boolean algebra. Indeed, we have 
a bit more: the center of any compact, ortho-coherent TOA -- and, in particular, any 
compact, lattice-ordered TOA -- is a compact Boolean algebra.

In [2], it is established that any TOML with an isolated point is discrete.
As Examples 2.4 and 2.5 illustrate, this is generally not true even 
for compact, lattice-ordered TOAs. On the other hand, a compact TOA in which $0$ is  isolated does have quite special properties, as we now show. 

Let us call an open set in a TOA {\em totally non-orthogonal} iff it contains no two orthogonal elements.

{\bf 2.10 Lemma:} {\em Every non-zero element of a TOA $L$ has a totally non-
orthogonal open neighborhood.}

{\bf Proof}: If $a \in L$ is non-zero, then $(a,a) \not \in \perp$. Since the
latter is closed in $L^{2}$, we can find open sets $U$ and $V$ with
$(a,a) \in U \times V$ and $(U \times V) \cap \perp = \emptyset$.
The set $U \cap V$ is a totally non-orthogonal open neighborhood of $a$.
$\Box$

{\bf 2.11 Proposition:} {\em Let $L$ be a compact TOA with $0$ isolated. Then
\begin{itemize}
\item[(a)] There is a finite bound on the size of pair-wise
orthogonal sets in $L$;
\item[(b)] Every chain in $L$ is finite;
\item[(c)] Every block of $L$ is finite;
\item[(d)] $L$ is atomic and of finite height;
\item[(e)] The set of atoms in $L$ is open.
\end{itemize}}

{\bf Proof}: \\
(a) If $0$ is isolated in $L$, then $L \setminus \{0\}$ is compact.
By Lemma 2.10, we can cover $L \setminus \{0\}$ by finitely many
totally non-orthogonal open sets $U_{1},...,U_{n}$. A pairwise-orthogonal
subset of $L \setminus \{0\}$ can meet each $U_{i}$ at most once, and so,
can have at most $n$ elements.

(b) Given an infinite chain in $L$, we can construct a
strictly increasing infinite sequence $x_{n}$ in $L$.
The sequence $y_{n}$ defined recursively by
$y_{1} = x_{1}$ and $y_{n} = x_{n+1} \ominus y_{n}$ for all $n \geq 2$, then gives us 
an infinite pairwise orthogonal set, contradicting part (a).

Now (c) follows immediately from (a), and (d) follows from (a) and (b) by considering
maximal chains. 

To establish (e), note that if $A$ and $B$ are any closed subsets
of $L$, then $(A \times B) \cap \perp$ is a closed, hence compact, subset of
$\perp$. Hence, the set $A \oplus B := \{ a \oplus b | a \in A, b \in B \ \text{and} \ a \perp b\}
= \oplus ((A \times B) \cap \perp)$ is closed. Now note that the set of non-atoms is
precisely $(L \setminus \{0\}) \oplus (L \setminus \{0\})$. Since $0$ is isolated,
$(L \setminus \{0\})$ is closed. Thus, the set of non-atoms is closed. $\Box$ 

If $a$ belongs to the center of $L$, then for each $b \in L$ there is a unique 
$c \in [0,a]$ with $(a \ominus c, c, b \ominus c) \in \M(L)$. This gives us a 
natural surjection $L \mapsto [0,a]$, which is easily seen to be continuous. Since 
the center of an orthoalgebra is a Boolean sub-orthoalgebra of $L$, and since a Boolean algebra of finite height is finite, Proposition 2.11 has the following 

{\bf 2.12 Corollary:} {\em Let $L$ be a compact TOA with $0$ isolated. Then the center of $L$ is finite. In particular, $L$ decomposes uniquely as the direct product of finitely many compact simple TOAs.}\\

\underline{\bf 3. Stably Ordered Topological Orthoalgebras}

In this section we consider a particularly tractable, but still quite broad, class of TOAs.

{\bf 3.1 Definition:} Call an ordered topological space $L$ {\em stably ordered} iff, for every open set $U \subseteq L$, the upper-set $U\uparrow = \{ b \in L | \exists a \in U a \leq b\}$ is again open.\footnote{The term used by Priestley [12] is ``space of type $I_{i}$."}

{\bf 3.2 Remarks:} Note that this is equivalent to saying that the second projection mapping $\pi_2 : \leq \rightarrow L$ is an open mapping, since for open sets $U, V \subseteq L$,
\[\pi_2(( U \times V) \cap \leq) = U\uparrow \cap V.\]
Note, too, that if $L$ carries a continuous orthocomplementation $'$, then $L$ is stably ordered iff $U\downarrow = \{ x | \exists y \in U, x \leq y\}$ be open for all open sets $U \subseteq L$.

{\bf 3.3 Example:} The following example (an instance of Example 2.3) shows that a TOA need not be stably ordered. Let $L = [0,1/4] \cup [3/4,1]$ with its usual topology, but without its usual order.
For $x, y \in L$, set $x \perp y$ iff $x + y = 1$ or $x = 0$ or $y = 0$. In any of these cases, define
$x \oplus y = x + y$. As is easily checked, this is a compact TOA. However, for the
clopen set $[0,1/4]$ we have $[0,1/4]\uparrow = [0,1/4] \cup \{1\}$, which is certainly not open.

The following is mentioned (without proof) in [12]:

{\bf 3.4 Lemma:} {\em Any topological $\wedge$-semilattice -- in particular, any topological lattice --  is stably ordered.}

{\bf Proof:} If $L$ is a topological meet-semilattice and $U \subseteq L$ is open, then
\[U\uparrow = \{ \ x \in L \ | \ \exists y \in U \ x \wedge y \in U\} = \pi_{1} (\wedge^{-1}(U))\]
where $\pi_{1} : L \times L \rightarrow L$ is the (open) projection map on the first factor and $\wedge : L \times L \rightarrow L$ is the (continuous) meet operation. $\Box$

Many orthoalgebras, including projection lattices, can be embedded in ordered abelian groups. Indeed, suppose $G$ is an ordered abelian group. If $e > 0$ in $G$, let $[0,e]$ denote the set of all elements $x \in G$ with $0 \leq x \leq e$. We can endow $[0,e]$ with the following partial-algebraic structure: for $x, y \in [0,e]$, set $x \perp y$ iff $x + y \leq e$, in which case let $x \oplus y = x + y$. Define $x' = e - x$. Then $([0,e],\oplus,',0,e)$ is an effect algebra (that is, satisfies all of the axioms for an orthoalgebra save possibly the condition that $x \perp x$ only for $x = 0$). By a {\em faithful sub-effect algebra} of $[0,e]$, we mean a subset $L$ of $[0,e]$, containing $0$ and $e$, that is closed under $\oplus$ (where this is defined) and under $'$, {\em and} such that, for all $x, y \in L$, $x \leq y$ iff $\exists z \in L$ with $y = x + z$. 

By way of example, let $L = L(\H)$, the projection lattice of a Hilbert space $\H$, regarded as an orthoalgebra, and let $G = {\cal B}_{sa}(\H)$, the ring of bounded self-adjoint operators on $\H$, ordered in the usual way. Then $L$ is a faithful sub-effect algebra of $[0,\1]$, where $\1$ is the identity operator on $\H$. This follows from the fact that, for projections $P, Q \in L(\H)$, $P + Q \leq \1$ iff $P \perp Q$, and the fact that if $P \leq Q$ as positive operators, then $Q - P$ is a projection. 
 
{\bf 3.5 Lemma:} {\em Let $L$ be an orthoalgebra, let $G$ be any ordered topological abelian group with a closed cone (equivalently, a closed order), and suppose that $L$ can be embedded as a sub-effect algebra of $[0,e]$, where $e > 0$ in $G$. Then $L$, in the 
topology inherited from $G$, is a stably ordered TOA.}

Proof: We may assume that $L$ is a subspace of $[0,e]$. Since $x \perp y$ in $L$ iff $x + y \leq e$, we have $\perp = +^{-1}([0,e]) \cap L$, which is relatively closed in $L$. The continuity of $\oplus$ and $'$ are automatic.
Suppose now that $U \cap L$ is a relatively open subset of $L$. Then, since $L$ is a faithful sub-effect algebra of $[0,e]$, the upper set generated by $U \cap L$ in $L$ is $U \uparrow \cap L$, where $U\uparrow$ is the upper set of $U$ in $[0,e]$. It suffices to show that this last is open. But $U \uparrow = \bigcup_{y \in G_{+}} U + y$, which is certainly open. $\Box$ 

It follows, in particular, that the projection lattice $L(\H)$ of a Hilbert space $\H$ is stably ordered in its norm topology. 

{\bf 3.6 Example:} A {\em state} on an orthoalgebra $(L,\oplus,0,1)$ is a mapping $f : L \rightarrow [0,1]$ such that $f(1) = 1$ and, for all $a, b \in L$,  $f(a \oplus b) = f(a) \oplus f(b)$ whenever $a \oplus b$ exists. A set $\Delta$ of states on $L$ is said to be {\em order-determining} iff $f(p) \leq f(q)$ for all $f \in \Delta$ implies $p \leq q$ in $L$. In this case the mapping $L \rightarrow {\Bbb R}^{\Delta}$ given by $p \mapsto \hat{p}$, $\hat{p}(f) = f(p)$, is an order-preserving injection. Taking $G = {\Bbb R}^{\Delta}$ in Lemma 3.5, we see that $L$ is a stably-ordered TOA in the topology inherited from pointwise convergence in $G$. As a special case, note that 
the projection lattice $L = L(\H)$ has an order-determining set of states of the 
form  form $f(p) = \langle p x, x \rangle$, where $x$ is a unit vector in $\H$.  Thus, $L(\H)$ is stably-ordered also in its weak topology.

If $U, V \subseteq L$, let us write $U \oplus V$ for $\oplus ((U \times V) \cap \perp)$, i.e., for the set of all (existing) orthogonal sums $a \oplus b$ with $a \in U$ and $b \in V$.

{\bf 3.7 Lemma:} {\em A TOA is stably ordered if, and only if, for every pair of open
sets $U, V \subseteq L$, the set $U \oplus V$ is also open.}

{\bf Proof:} Let $U$ and $V$ be any two open sets in $L$. Then
\begin{eqnarray*}
U \oplus V & = & \{ c \in L | c = a \oplus b, a \in U, b \in V\}\\
& = & \{ c \in L | \exists a \in U \ a \leq c \ \text{and} \ c \ominus a \in V\} \\
& = & \pi_{2}(\ominus^{-1}(V) \cap (L \times U \uparrow)) .
\end{eqnarray*}
Now, since $L$ is stably ordered, $U\uparrow$ is open, and hence, $\ominus^{-1}(V) \cap (L \times U \uparrow)$ is
relatively open in $\leq$. But as observed above, for $L$ stably ordered, $\pi_2 : \leq \rightarrow L$ is an open mapping,
so $U \oplus V$ is open. For the converse, just note that $U \uparrow = U \oplus L$. $\Box$

Proposition 2.11 tells us that a compact TOA $L$ with $0$ isolated is atomic and of finite height. It follows easily that every element of $L$ can be expressed as a finite orthogonal sum of atoms. Let the {\em dimension}, $\dim(a)$, of an element $a \in L$ be the minimum number $n$ of atoms $x_1,...,x_n$ such that $a = x_1 \oplus \cdots \oplus x_n$. Note also that $a \in L$ is an atom iff $\dim(a) = 1$. 

{\bf 3.8 Theorem:} {\em Let $L$ be a compact, stably-ordered TOA in which $0$ is an isolated point.
Then
\begin{itemize}
\item[(a)] The dimension function is continuous -- equivalently, the set of elements of a given dimension is clopen.
\item[(b)] The topology on $L$ is completely determined by that on the set of atoms. \end{itemize}}

{\bf Proof:} By Lemma 3.7 and the fact that $L$ is compact, if $A$ and $B$ are clopen subsets of $L$, then $A \oplus B$ is again clopen. Since $0$ is isolated, $L \setminus \{0\}$ is clopen. Since the set of non-atoms in $L$ is exactly $(L \setminus \{0\}) \oplus (L \setminus \{0\})$, it follows that the set of atoms is clopen. Now define
a sequence of sets $L_{k}$, $k = 0,...,\dim(1)$, by setting by $L_{0} = \{0\}$, $L_{1}$ = the set of atoms of $L$,
and $L_{k+1}:= L_{k} \oplus L_{1}$. These sets are clopen, as are all Boolean combinations of them.
In particular, the set
\[\{ a \in L | \dim(a) = k\}  = L_{k} \setminus ( \bigcup_{i=0}^{k-1} L_{k})\]
is clopen for every $k = 0,...,\dim(1)$. This proves (a). For (b), it now suffices to show that the topology on each $L_n$, $n > 1$, is determined by that on $L_1$. Since $L_1$ and $L_2$ are clopen, Lemma 3.7 tells us that the mapping $\oplus : ((L_n \times L_1) \cap \perp) \rightarrow L_{n_+1}$ is an open surjection, and hence, a quotient mapping. Thus, the topology on $L_{n+1}$ is entirely determined by that on $L_{n}$ and that on $L_{1}$. An easy induction completes the proof. $\Box$ \\

\underline{\bf 4. The Logic of an Algebraic Topological Test Space}

We now take up the problem of topologizing the representation
of orthoalgebras as ``logics" of algebraic test space, glossed in section 1. 
This continues work begun in [16], which  we briefly reprise here.

{\bf 4.1 Definition:} A {\bf topological
test space} is a test space $(X,{\frak A})$ where \begin{itemize} \item[(a)]
$X$ is a Hausdorff space, \item[(b)] The relation $\perp \subseteq X^{2}$
is closed. \end{itemize}

Condition (b) guarantees that every outcome $x \in X$ has an local basis of
neighborhoods $V$ that are totally non-orthogonal, in that $(V \times
V) \cap \perp = \emptyset$. Thus, pairwise orthogonal sets -- and
in particular, tests -- are discrete. It can also be shown ([16], Lemma 2.3) that they are closed.

If $X$ is any topological space and ${\frak A}$ is any collection of closed subsets
of $X$, the {\em Vietoris topology} on $\frak A$ is
that generated by sets of the form
\[[U] := \{ E \in {\frak A} | E \cap U \not =
\emptyset\}\ \text{and} \ (U) := [U^{c}]^{c} = \{ E \in {\frak A} | E
\subseteq U\}\]
where $U \subseteq X$ is open.
We denote by $2^{X}$ the hyperspace of all closed subsets of $X$. It is a basic fact [10] that $X$ is compact iff $2^{X}$
is compact. \footnote{Also, if $X$ is a compact metric space, the Vietoris topology on
$2^{X}$  coincides with
the Hausdorff metric topology.}

If $(X,{\frak A})$ is any topological test space,
we regard both ${\frak A}$ and ${\cal E}(X,{\frak A})$ as subspaces of $2^{X}$.  We can
endow the logic $\Pi(X,{\frak A})$ with the quotient topology induced by the
the natural surjection $p : {\cal E}(X,{\frak A}) \rightarrow \Pi(X,{\frak A})$.
In [16], we obtained a sufficient condition for $\Pi(X,{\frak A})$, with this quotient topology, to be a topological orthoalgebra. We now ask: given a topological OA $L$, when is the natural orthoalgebra isomorphism $\Pi({\frak A}_L) \rightarrow L$ actually an homeomorphism? Equivalently, when is the canonical surjection $q : {\cal E}(X_L,{\frak A}_L) \rightarrow L$ given by $q: A \mapsto \bigoplus A$, a quotient mapping? We shall show (Theorem 4.6 below) that this is always the case if $L$ is stably ordered. 

For the remainder of this section, $L$ denotes a topological orthoalgebra, $(X_L, {\frak A}_L)$, the associated test space of finite orthopartitions of the unit in $L$, and ${\cal E}$, the set of events for the latter, i.e., the set of all finite summable subsets
of $L$.

{\bf 4.2 Lemma:} {\em Let $L$ be any topological orthoalgebra, and let ${\cal E}^{(n)}$ denote the set of $n$-element
events in ${\cal E}$. Then
\begin{itemize}
\item[(a)] ${\cal E}^{(n)}$ is clopen in ${\cal E}$
\item[(b)] The relative Vietoris topology on ${\cal E}^{(n)}$ is generated by sets of the form $[U]$, $U \subseteq L$ open.
\end{itemize}}

{\bf Proof:} (a) Let $A = \{a_1,...,a_n\} \in {\cal E}^{(n)}$. For each $a_i$ in $A$, let $U_i$ be a totally non-orthogonal open neighborhood. Since $L$ is Hausdorff, we can assume that these are pairwise disjoint. Now any event in
the basic Vietoris open set $\langle U_1,...,U_n \rangle$ has exactly $n$ elements, so ${\cal U} \subseteq {\cal E}^{(n)}$.
This shows the latter is open. Since the sets ${\cal E}^{(n)}$ are pairwise disjoint and cover ${\cal E}$, each is clopen.
For (b), note that $\langle U_1,...,U_n \rangle =  [U_1] \cap \cdots \cap [U_n] \cap {\cal E}^{(n)}$. $\Box$

{\bf 4.3 Lemma:} {\em The canonical surjection $q : {\cal E} \rightarrow L$ is continuous.}

{\bf Proof:} In view of Lemma 4.2, it suffices to show that the restriction of $q$ to each clopen set ${\cal E}^{(n)}$ is continuous. Let ${\cal E}^{(n)}_{o} \subseteq L^{n}$ denote the set of $n$-tuples $(a_1,...,a_n)$ for which $a_1 \oplus \cdots \oplus a_n$ exists. There is a natural quotient mapping $\pi : {\cal E}^{(n)}_{o} \rightarrow {\cal E}^{(n)}$ that forgets the order. It is not hard to show [14] that this is an open continuous mapping for the relative product topology on ${\cal E}^{(n)}_{o}$ and the Vietoris topology on ${\cal E}^{(n)}$. Note that ${\cal E}^{(2)}_{o} = \perp$. Using the continuity of $\oplus$ on $\perp$, plus the strong associativity of $\oplus$ as a partial binary operation, one can show by induction that the mapping $\oplus^{(n)} : {\cal E}^{(n)}_{o} \rightarrow L$ given by $\oplus^{(n)} : (a_1,....,a_n) \mapsto a_1 \oplus \cdots \oplus a_n$ is continuous on ${\cal E}^{(n)}_{o}$. Since $\oplus^{(n)} = q \circ \pi$, it follows that $q$ is continuous. Indeed, for any open set $U \subseteq L$, ${\oplus^{(n)}}^{-1}(U) = \pi^{-1}(q^{-1}(U))$ is open. Since $\pi$ is open and surjective, therefore, $q^{-1}(U) = \pi(\pi^{-1}(q^{-1}(U)))$ is also open. $\Box$.

We are now in a position to prove the advertised result:

\newpage
{\bf 4.5 Theorem:} {\em If $L$ is stably ordered, then
the canonical surjection $q : {\cal E} \rightarrow L$ is open, and hence, a quotient mapping. }

{\bf Proof:} By part (b) of Lemma 4.2, a basis for the Vietoris topology on ${\cal E}^{(n)}$ consists of sets of the form
$[U_1] \cap \cdots \cap [U_n] \cap {\cal E}^{(n)}$,
where $U_1,...,U_n$ are totally non-orthogonal open sets in $L$. Now, an event $\{a_1,...,a_n\}$ belonging to this set is
just a selection of $n$ jointly orthogonal elements $a_1 \in U_1, ....,a_n \in U_n$, so
 $q\left ( [U_1] \cap \cdots \cap [U_n] \cap {\cal E}^{(n)} \right ) \ = \ U_1 \oplus \cdots \oplus U_n$.
By Lemma 3.7, this last set is open in $L$. Thus, the restriction of $q$ to  ${\cal E}^{(n)}$ is an open mapping. By part (a) of Lemma 4.2, ${\cal E}^{(n)}$ is clopen in $L$.  Since ${\cal E} = \bigcup_{n} {\cal E}^{(n)}$, the mapping $q : {\cal E} \rightarrow L$ is open as well. By Lemma 4.3, $q$ is continuous, and hence, a quotient mapping. In particular, the quotient topology on $L$ agrees with the original topology.  $\Box$

{\bf Remark:} The converse also holds, since if $q$ is open, then so is $q([U]) = U\uparrow$ for any open set $U \subseteq L$.

{\bf 4.6 Corollary:} {\em Any stably-ordered topological orthoalgebra $L$ is canonically isomorphic, both algebraically and topologically, to the logic $\Pi(X_L,{\frak A}_L)$ of 
its associated space of orthopartitions.}

In [16], we called a topological test space $(X,{\frak A})$ {\em
stably complemented} if the set ${\cal U}^{\oc} = \{ B \in {\cal
E} | \exists A \in {\cal U} , A \ \text{complementary to} \ B\}$ is open for every
Vietoris-open set ${\cal U}$ of events. We showed ([16], Lemma 3.5
and Proposition 3.6) that $(X,{\frak A})$ is stably complemented
iff $p : {\cal E} \rightarrow \Pi$ is open and $' : \Pi
\rightarrow \Pi$ is continuous, and that if $(X,{\frak A})$ is
stably complemented with ${\cal E}$ closed in $2^X$, then $\Pi$ is a
topological orthoalgebra. 

{\bf 4.7 Proposition:} {\em Let $L$ be a topological orthoalgebra. Then
\begin{itemize}
\item[(a)] If $L$ is stably ordered, then $(X_L, {\frak A}_L)$ is stably complemented.
\item[(b)] If $L$ is compact with $0$ isolated,
and $(X_L, {\frak A}_L)$ is stably complemented, then $L$ is stably ordered. \end{itemize}}

To prove this, we require a preliminary 

{\bf 4.8 Lemma:} {\em For each $n$, the set ${\cal E}^{(n)}$ of $n$-element events is closed in $2^{L}$ in the Vietoris topology.}

{\bf Proof:} Note first that if $A \in {\cal E}$ and $x \perp \bigoplus A$, then $A \cup \{x\}$ is again an event. We now proceed by induction, noting that the result is trivial for $n = 1$. Supposing it to hold for a given $n$, consider a net of events $A_{\lambda}$ in ${\cal E}^{(n+1)}$ converging (in the Vietoris topology) to a set $A \in 2^{L}$. For any $k$ distinct elements $x_1,...,x_k$ of $A$, select pairwise disjoint open neighborhoods $U_1,...,U_k$. Then $A \in [U_1] \cap \cdots \cap [U_k]$; hence, 
$A_{\lambda}$ must eventually meet each of the disjoint sets $U_1,...,U_k$. It 
follows that $k$ can be no larger than $|A_{\lambda}| = n+1$. Hence, $A$ itself is a finite set with $|A| \leq n + 1$. But if we now assume that each $U_{i}$ is totally non-orthogonal, we see that $A_{\lambda}$ can meet each $U_i$ at most once, whence, $|A| =  n+1$. To see that $A$ is in fact an event, note that (eventually) we have a canonical
bijection $\phi_{\lambda} : A \mapsto A_{\lambda}$  sending each $x \in A$ to the unique element $x_\lambda$ of $A_{\lambda} \cap U$, $U$ any totally non-orthogonal neighborhood of $x$ containing no other point of $A$. Moreover, the net $x_{\lambda}$ converges to $x$. Now select any $x \in A$, and let $B = A \setminus x$. Let $x_{\lambda}$ be as above, and let $B_{\lambda} = \phi_{\lambda}(B)$. It is easy to check that $B_{\lambda} \rightarrow B$ in the Vietoris topology.
By the continuity of $\bigoplus$ (Lemma 4.3), we have $\bigoplus B_{\lambda} \rightarrow \bigoplus B$. Since $x_{\lambda} \perp \bigoplus B_{\lambda}$ and the relation $\perp$ is closed, we have $x \perp \bigoplus B$ -- whence, $B \cup \{x\} = A$ is an event. $\Box$

{\bf Proof of Theorem 4.7:} (a) If $L$ is stably ordered, then by Lemma 4.2 and Theorem 4.3, the canonical mapping $q : {\cal E} \rightarrow L$ is both continuous and open. Hence,
for any Vietoris-open set ${\cal U}$, we have $q({\cal U})$ open. Since the orthocomplementation $' : L \rightarrow L$ is
an homeomorphism, $q({\cal U})'$ is also open, whence, so is $q^{-1}(q({\cal U})') = {\cal U}^{\oc}$.

(b) Let $\phi : \Pi \rightarrow L$ be the (well-defined) bijection
sending $p(A)$ to $\bigoplus A = q(A)$ for all $A \in {\cal E}$.
We have the following commutative diagram, with $p$ open, $q$
continuous and $\phi$, a bijection:
\[\begin{array}{ccc}
 & {\cal E} & \\
\ \ \ \ \ \ \ \ \ \ \ \ \ \ \  \ \ \ \ \  \ \ \ \ \ \ \ \ \ \ \ ^{q}\swarrow &  &  \searrow ^{p} \ \ \ \ \ \ \ \ \ \ \ \ \ \ \  \ \ \ \ \  \ \ \ \ \ \ \ \ \ \  \\
\ \ \ \ \ \ \ \ \ \ \ \ \ \ \  \ \ \ \ \ L & \stackrel{\phi}{\longleftarrow} & \Pi \ \ \ \ \ \ \ \ \ \ \ \ \ \ \ \ \ \  \ \ \end{array} \]
Hence, for any $U \subseteq L$ open, we have $p^{-1}\phi^{-1}(U) = (\phi \circ p)^{-1} (U) = q^{-1}(U)$, which
is open. But then $\phi^{-1}(U)$ is open, since $p$ is a quotient mapping. So $\phi$ is continuous. Since $L$ is
compact, $2^{L}$ is also compact. Since $L$ is compact with $0$ is isolated in $L$, it 
follows from Proposition 2.10 and Lemma 4.8 that ${\cal E}$ is the union of finitely many closed -- hence, compact -- sets ${\cal E}^{(n)}$. Thus, $\cal E$ is itself compact, whence, so is $\Pi$. It follows that the continuous bijection $\phi$ is in fact a homeomorphism. $\Box$\\

{\bf References}

[1] Bennett, M. K., and Foulis, D. J., {\em Interval effect algebras and unsharp quantum logics}, Advances in Applied Math. {\bf 19} (1997), 200-219.

[2] Choe, T. H., and Greechie, R. J., {\em Profinite Orthomodular Lattices},
Proc. Amer. Math. Soc. {\bf 118} (1993), 1053-1060

[3] Choe, T. H., Greechie, R. J., and Chae, Y., {\em Representations of locally compact
orthomodular lattices}, Topology and its Applications {\bf 56} (1994) 165-
173

[4] Foulis, D. J., Greechie, R. J.,
and Ruttimann, G. T., {\em Filters and Supports on  Orthoalgebras} Int. J.
Theor. Phys. {\bf 31} (1992) 789-807

[5] Foulis, D. J., Greechie, R. J., and
R\"{u}ttimann, G. T., {\em Logico-algebraic structures  II: supports
on test spaces}, Int. J. Theor. Phys.
{\bf 32} (1993) 1675-1690.

[6] Foulis, D. J., and Randall, C.H., {\em What are quantum logics, and what ought they to be?}, in Beltrametti, E., and van Fraassen, B. C. (eds.), {\bf Current Issues in Quantum Logic}, Plenum: New York, 1981. )

[7] Greechie, R. J., Foulis, D. J., and Pulmannov\'{a}, S., {\em The center
of an effect algebra}, {\em Order} {\bf 12} (1995), 91-106.

[8] Johnstone, P. T., {\bf Stone Spaces}, Cambridge: Cambridge University Press, 1982

[9] Kalmbach, G., {\bf Orthomodular Lattices}, Academic Press, 1983

[10] Michael, E., {\em Topologies on spaces of subsets},
Trans. Am. Math. Soc. {\bf 7} (1951) 152-182

[11] Nachbin, L., {\bf Topology and Order}, van Nostrand: Princeton 1965

[12] Priestley, H. A.,  {\em Ordered Topological Spaces and the Representation of Distributive Lattices}, Proc. London Math. Soc. {\bf 24} (1972), 507-530. 

[13] Pulmannov\'{a}, S., and Riecanova, Z., {\em Block-finite orthomodular
lattices}, J. Pure and Applied Algebra {\bf 89} (1993), 295-304

[14] Pulmannov\'{a}, S., and Rogalewicz, V., {\em Orthomodular lattices with almost
orthogonal sets of atoms}, Comment. Math.Univ. Carolinae {\bf 32} (1991)
423-429.

[15] Wilce, A., {\em Test Spaces and Orthoalgebras}, in Coecke et al (eds.)
{\bf Current Research in Operational  Quantum Logic}, Kluwer: Dordrecht
(2000).

[16] Wilce, A., {\em Topological Test Spaces}, to appear, International Journal of Theoretical Physics.

\end{document}